# Correction to: The scaling limit behaviour of periodic stable-like processes

BY BRICE FRANKE. *Bernoulli*: (2006), **12**, 551–570

An error occurs in the proof of Proposition 3. The term $h(n, y)$ introduced in the last line of page 558 does not converge to zero after division by $g(n)$.

For Proposition 3 to hold, the distribution function $F_\alpha := \pi \circ \alpha^{-1}$ must satisfy the following additional assumption:

*The distribution function $F_\alpha$ has a density $F'_\alpha$ and there exists $\beta > 0$ such that*

$$\frac{F'_\alpha(tp + \alpha_o)}{F'_\alpha(t + \alpha_o)} \longrightarrow p^\beta \qquad as\ t \to 0.$$

The co-area formula from geometric measure theory and the fact that $\pi$ has a density with respect to the Lebesgue measure ensure the existence of $F'_\alpha$ in our situation. The exponent $\beta$ is related to the dimension of the set where $\alpha$ is equal to $\alpha_o$. Proposition 3 holds under the previous assumption if the scale function $g(n)$ is replaced by the new scale function

$$\tilde{g}(n) := F'_\alpha\left(\frac{\alpha_o}{\log n} + \alpha_o\right) \frac{\alpha_o}{\log n} \Gamma(\beta + 1).$$

In the statements of Corollary 1, Theorem 3 and Theorem 4, the scale function $g(n)$ has also to be replaced by $\tilde{g}(n)$.

The following modification has to be made in the proof of Proposition 3:

For $k_n := \alpha_o^{-1}(2 - \alpha_o) \log n$ and $m_n := \alpha_o / \log n$, we have

$$\int_{\alpha_o}^{2} |y|^{-1-q} n^{1-q/\alpha_o} F_\alpha(\mathrm{d}q) = \int_0^{2-\alpha_o} |y|^{-1-r-\alpha_o} n^{-r/\alpha_o} F'_\alpha(r + \alpha_o) \,\mathrm{d}r$$

$$= m_n \int_0^{k_n} |y|^{-1-m_n p - \alpha_o} e^{-p} F'_\alpha(m_n p + \alpha_o) \,\mathrm{d}p$$

$$= m_n F'_\alpha(m_n + \alpha_o) |y|^{-1-\alpha_o} \Gamma(\beta + 1) + E(n, y)$$

$$= g(n)|y|^{-1-\alpha_o} + E(n, y).$$

The assumption on $F'_\alpha$ implies that $g(n)^{-1}E(n,y) \to 0$ as $n \to \infty$. We note that the convergence is uniform in $y$ on the exterior of an arbitrary small ball centered at zero. From this point, the proof given in the article remains unchanged.